\documentclass[11pt]{article}

\usepackage[T1]{fontenc}
\usepackage{lmodern}
\usepackage{microtype}
\usepackage{amsmath,amssymb,amsthm,mathtools}
\usepackage{enumitem}
\usepackage{geometry}
\usepackage{hyperref}
\usepackage[nameinlink,noabbrev]{cleveref}
\usepackage{booktabs}
\usepackage{array}
\usepackage{xcolor}
\usepackage{bm}

\hypersetup{colorlinks=true,linkcolor=blue,citecolor=blue,urlcolor=blue}
\setlist[itemize]{leftmargin=2em}
\setlist[enumerate]{leftmargin=2.2em}

\newtheorem{theorem}{Theorem}[section]
\newtheorem{lemma}[theorem]{Lemma}
\newtheorem{proposition}[theorem]{Proposition}
\newtheorem{corollary}[theorem]{Corollary}
\crefname{theorem}{Theorem}{Theorems}
\crefname{lemma}{Lemma}{Lemmas}
\crefname{proposition}{Proposition}{Propositions}
\crefname{corollary}{Corollary}{Corollaries}
\theoremstyle{definition}

\newtheorem{example}[theorem]{Example}
\theoremstyle{remark}
\newtheorem{remark}[theorem]{Remark}

\newcommand{\Leb}{\operatorname{Leb}}
\newcommand{\TV}{\operatorname{TV}}
\newcommand{\BL}{\operatorname{BL}}

\newcommand{\ind}{\mathbf{1}}

\newcommand{\Mhat}{\widehat M}
\newcommand{\muHat}{\widehat\mu}
\newcommand{\dd}{\,\mathrm{d}}
\newcommand{\norm}[1]{\left\lVert #1\right\rVert}

\title{Statistical Stability and Linear Response for Suspension Measures under Base and Roof Perturbations}
\author{
Aparna Rajput\\
Department of Mathematics and Statistics, Concordia University\\
\texttt{aparna.ar.rajput@gmail.com}
}
\date{November 2, 2026}

\begin{document}
\maketitle

\begin{abstract}
We study the dependence of canonical invariant probability
measures of suspension semiflows on simultaneous perturbations
of the invariant base density and the roof function. We compare
two natural representations of the suspension measure:
physical-height coordinates, in which the support depends on
the roof function, and normalised-height coordinates, in which
all measures are defined on a fixed cylinder.

In physical-height coordinates, we obtain an exact fibrewise
formula for the $L^1$ discrepancy between suspension densities
and derive quantitative total-variation estimates. In
normalised-height coordinates, we establish corresponding
stability estimates through the product of the base density
and the roof function. These results transfer quantitative
statistical-stability estimates for the base densities to the
associated suspension measures.

For one-parameter families whose base densities are
differentiable in $L^1$ and whose roof functions are
differentiable in $L^\infty$, we derive an $L^1$
linear-response formula on the fixed cylinder and a
bounded-Lipschitz first-variation formula in physical
height.
The latter contains a signed measure supported on the graph
of the unperturbed roof. We further compute the exact
first-order total-variation speed and prove that
total-variation differentiability holds if and only if
the first-order roof velocity vanishes almost everywhere
with respect to the unperturbed fibre density. The results
are measure-theoretic and do not require spectral or mixing
assumptions on the suspension semiflow.
\end{abstract}



\section{Introduction}

Suspension constructions convert discrete-time dynamics into continuous-time
dynamics by assigning a positive return time, or roof function, to each point
of the base. If a measurable map $T$ preserves a probability measure $\nu$
and the roof $\tau$ is positive and $\nu$-integrable, the canonical suspension
measure is obtained by taking the product of $\nu$ with Lebesgue measure
along fibres and normalising by the mean roof height. This elementary formula
is standard, but perturbing both $T$ and $\tau$ raises a comparison issue that
is easy to overlook: changing the roof changes not only the density but also
the physical suspension domain.

The purpose of this paper is to isolate the measure-theoretic consequences
of this moving-domain geometry and to separate them from the substantially
harder spectral theory of suspension semiflows. For piecewise-expanding
maps, transfer-operator methods and Lasota--Yorke inequalities provide a
classical route to absolutely continuous invariant measures and their
stability; see, for example,
\cite{LasotaYorke1973,Keller1985,KellerLiverani1999,RajputGora2025,GoraRajput2026}.
Quantitative perturbation estimates for piecewise-expanding families have
been developed in several settings, including families with moving branch
geometry \cite{GalatoloLucena2020,AlvesEtubi2025}. By contrast, spectral and
correlation results for suspension semiflows require additional structure
in the flow direction and non-arithmetic or transversality hypotheses on
the roof; see \cite{Tsujii2008,DaltroVarandas2021,Butterley2014}.

The contribution of the present work is different from the
spectral and mixing theory of suspension semiflows. We
investigate the dependence of the canonical suspension
probability measure on simultaneous perturbations of the
invariant base density and the roof function. In particular,
we distinguish between physical-height coordinates, where
motion of the roof boundary produces a singular contribution
to the first variation, and normalised-height coordinates,
where this moving-boundary contribution is absent. This
distinction leads to different differentiability properties
of the suspension measures. We obtain an exact formula for
the first-order total-variation speed and a necessary and
sufficient condition for differentiability in the
total-variation norm. Our results also provide a framework
for transferring quantitative stability and linear-response
information from the base density to the associated
suspension measure.

Our results concern the invariant suspension measure itself. They are valid
for a general measurable base once its invariant density is given. The
dynamical assumptions needed to produce or perturb that density enter only
later, through applications.

There are two natural comparison geometries.
\begin{itemize}
\item In \emph{physical height}, a suspension is represented in the ambient
strip $I\times[0,\infty)$ by the subgraph $0\le s<\tau(x)$. Perturbing
$\tau$ moves the support boundary.
\item In \emph{normalised height}, the coordinate $u=s/\tau(x)$ identifies
every suspension with the fixed cylinder $\Mhat=I\times[0,1)$. The Jacobian
of this change of variables inserts the factor $\tau(x)$ into the transported
density.
\end{itemize}

These geometries lead to different notions of response. On the
fixed cylinder, if the base densities are differentiable in
$L^1$ and the roof functions are differentiable in $L^\infty$,
then the transported suspension densities have an ordinary
$L^1$ derivative. In physical height, the moving boundary
contributes a graph-supported signed measure to the first
variation. Thus bounded-Lipschitz response can coexist with
failure of differentiability in the total-variation norm.
The differentiability of the base density itself is not
established here in general; rather, it is assumed when
deriving the corresponding response formulas for the
suspension measure.

The main results are measure-theoretic. We first establish
an exact formula for the $L^1$ discrepancy between two
physical suspension densities and derive quantitative
total-variation estimates. We then show that, after fibre
normalisation, the comparison reduces to the difference
between the products of the base densities and roof
functions. For differentiable families, we obtain
first-variation formulas in both representations. The
main distinction appears in physical-height coordinates:
the first variation contains a singular measure supported
on the graph of the unperturbed roof. We use this term to
compute the exact first-order total-variation speed and
to characterise precisely when the physical suspension
measure is differentiable in total variation. If
$k_0=h_0/A_0$ denotes the physical fibre density and $\dot\tau$ the
first-order roof velocity, then
\[
 \lim_{\theta\to0}
 \frac{d_{\TV}(\mu_\theta^{\rm ph},\mu_0^{\rm ph})}{|\theta|}
 =
 \frac12\left(
 \int_I \tau_0|\dot k|\,\dd x+
 \int_I k_0|\dot\tau|\,\dd x
 \right),
\]
and the curve is differentiable in total variation if and only if
$k_0\dot\tau=0$ almost everywhere. In particular, whenever $h_0>0$ almost
everywhere, a genuinely moving roof prevents total-variation linear response
even though the curve remains locally Lipschitz in total variation.

\section{Suspension measures and invariance}

Let $I=[0,1]$ with Lebesgue measure. Let $T:I\to I$ be measurable and
suppose that $T$ preserves an absolutely continuous probability measure
\[
 \dd\nu(x)=h(x)\,\dd x,\qquad
 h\ge0,\qquad
 \int_I h\,\dd x=1.
\]
Let $\tau:I\to(0,\infty)$ be measurable and suppose
\begin{equation}\label{eq:roof-lower}
 \tau(x)\ge m>0
 \qquad\text{for Lebesgue-a.e. }x.
\end{equation}
Assume also
\begin{equation}\label{eq:meanroof}
 0<A:=\int_I h(x)\tau(x)\,\dd x<\infty.
\end{equation}
The upper boundedness of $\tau$ is not required in this section.

Because \cref{eq:roof-lower} is only an almost-everywhere statement, it is
convenient to remove once and for all the null set that could be reached by
a forward orbit.

\begin{lemma}\label{lem:safe}
Let
\[
 N=\{x\in I:\tau(x)<m\}
\]
and define
\[
 X_*:=I\setminus\bigcup_{j\ge0}T^{-j}N.
\]
Then $\nu(X_*)=1$, $T(X_*)\subset X_*$, and
\[
 \tau(T^j x)\ge m
 \qquad\text{for every }x\in X_*\text{ and }j\ge0.
\]
\end{lemma}

\begin{proof}
Since $\nu\ll\Leb$ and $\Leb(N)=0$, we have $\nu(N)=0$. By
$T$-invariance of $\nu$,
\[
 \nu(T^{-j}N)=\nu(N)=0
 \qquad(j\ge0).
\]
A countable union of $\nu$-null sets is null, hence $\nu(X_*)=1$.
If $x\in X_*$ and $Tx\notin X_*$, then
$T^{j+1}x\in N$ for some $j\ge0$, contradicting $x\in X_*$. The final
assertion follows directly from the definition of $X_*$.
\end{proof}

Define the measurable suspension domain
\[
 M_{\tau,*}
 :=
 \{(x,s):x\in X_*,\ 0\le s<\tau(x)\},
\]
with the usual boundary identification
\[
 (x,\tau(x))\sim(Tx,0).
\]
Put $\tau_0(x)=0$ and, for $n\ge1$,
\[
 \tau_n(x)
 =
 \sum_{j=0}^{n-1}\tau(T^j x).
\]
By \cref{lem:safe},
\[
 \tau_n(x)\ge nm
\]
on $X_*$, so $\tau_n(x)\to\infty$. Consequently, for each
$(x,s)\in M_{\tau,*}$ and $t\ge0$ there is a unique
$n=n(x,s,t)\ge0$ such that
\[
 \tau_n(x)\le s+t<\tau_{n+1}(x).
\]
Define
\begin{equation}\label{eq:semiflow}
 \Phi^t(x,s)
 =
 \bigl(T^n x,\ s+t-\tau_n(x)\bigr).
\end{equation}
This gives the usual forward suspension semiflow on a full-$\nu$-measure
base set. If desired, one may instead replace $\tau$ on a null set by a
representative satisfying the lower bound everywhere.

The canonical probability measure is
\begin{equation}\label{eq:canonical}
 \mu(B)
 =
 \frac1A
 \int_I\int_0^{\tau(x)}
 \ind_B(x,s)\,\dd s\,\dd\nu(x)
 =
 \frac1A
 \int_I\int_0^{\tau(x)}
 \ind_B(x,s)h(x)\,\dd s\,\dd x.
\end{equation}
Thus, on the physical suspension domain, the density with respect to
$\dd x\,\dd s$ is $h(x)/A$.

\begin{proposition}[Invariance]\label{prop:invariance}
The probability measure $\mu$ in \cref{eq:canonical} is invariant under
$(\Phi^t)_{t\ge0}$.
\end{proposition}

\begin{proof}
It suffices to prove invariance for $0\le t<m$ and then use the semigroup
property. Let $g$ be bounded and measurable. For $t<m$, an orbit can cross
at most one roof during a time interval of length $t$: after a crossing the
new fibre has height at least $m>t$.

Splitting each fibre at $\tau(x)-t$ and changing the vertical variable gives
\begin{align*}
 &\int_{X_*}\int_0^{\tau(x)}
 g(\Phi^t(x,s))\,\dd s\,\dd\nu(x)
 \\
 &=
 \int_{X_*}\int_t^{\tau(x)}
 g(x,r)\,\dd r\,\dd\nu(x)
 +
 \int_{X_*}\int_0^t
 g(Tx,r)\,\dd r\,\dd\nu(x).
\end{align*}
For each fixed $r$, $T$-invariance of $\nu$ and forward invariance of $X_*$
imply
\[
 \int_{X_*}g(Tx,r)\,\dd\nu(x)
 =
 \int_{X_*}g(x,r)\,\dd\nu(x).
\]
Hence the preceding sum equals
\[
 \int_{X_*}\int_0^{\tau(x)}
 g(x,r)\,\dd r\,\dd\nu(x).
\]
After division by $A$, invariance holds for $0\le t<m$. Every $t\ge0$ is a
finite sum of positive times smaller than $m$, so invariance follows from
the semiflow property.
\end{proof}

\begin{remark}\label{rem:integrability}
The essential integrability condition for the canonical probability
measure is precisely
\[
 h\tau\in L^1(I),
 \qquad
 0<\int_Ih\tau\,\dd x<\infty.
\]
Integrability of $\tau$ with respect to Lebesgue measure alone is neither
necessary nor, for an arbitrary density $h$, sufficient.
\end{remark}

\section{Stability estimates for suspension measures}

We first record a normalisation inequality used repeatedly below. For a
finite signed measure $\sigma$, write
\[
 \norm{\sigma}_{\mathrm{var}}
 :=
 |\sigma|(X)
\]
for its total variation norm. For probability measures $\rho$ and $\eta$
on the same measurable space, we use
\[
 d_{\TV}(\rho,\eta)
 :=
 \sup_B|\rho(B)-\eta(B)|
 =
 \frac12\norm{\rho-\eta}_{\mathrm{var}}.
\]

\begin{lemma}[Normalisation inequality]\label{lem:normalization}
Let $(X,\mathcal A,m)$ be a measure space and let
$f_0,f_1\in L^1(m)$ be nonnegative with masses
\[
 A_i=\int_X f_i\,\dd m>0,
 \qquad
 \rho_i=\frac{f_i}{A_i}m.
\]
Then
\begin{equation}\label{eq:norm-main}
 d_{\TV}(\rho_1,\rho_0)
 \le
 \frac{\norm{f_1-f_0}_1+|A_1-A_0|}
 {2\max\{A_0,A_1\}}
 \le
 \frac{\norm{f_1-f_0}_1}
 {\max\{A_0,A_1\}}.
\end{equation}
If $A_0=A_1=A$, then the exact identity
\begin{equation}\label{eq:norm-equal-mass}
 d_{\TV}(\rho_1,\rho_0)
 =
 \frac1{2A}\norm{f_1-f_0}_1
\end{equation}
holds.
\end{lemma}

\begin{proof}
Assume without loss of generality that $A_1\ge A_0$. Adding and subtracting
$f_0/A_1$ gives
\begin{align*}
 \norm{\frac{f_1}{A_1}-\frac{f_0}{A_0}}_1
 &\le
 \frac1{A_1}\norm{f_1-f_0}_1
 +
 \left|
 \frac1{A_1}-\frac1{A_0}
 \right|\norm{f_0}_1
 \\
 &=
 \frac{
 \norm{f_1-f_0}_1+|A_1-A_0|
 }{A_1}.
\end{align*}
Dividing by $2$ proves the first inequality. Since
\[
 |A_1-A_0|
 \le
 \norm{f_1-f_0}_1,
\]
the second follows. If $A_0=A_1=A$, the two normalised densities have the
same denominator, yielding \cref{eq:norm-equal-mass} directly.
\end{proof}

\subsection{Physical-height comparison}

Let $(h_i,\tau_i)$, $i=0,1$, be two pairs with
\[
 h_i\ge0,
 \qquad
 \int_Ih_i\,\dd x=1,
 \qquad
 \tau_i>0,
\]
and
\[
 0<A_i
 :=
 \int_I h_i\tau_i\,\dd x
 <\infty.
\]
Embed both suspensions in the common ambient half-strip
\[
 \mathcal R:=I\times[0,\infty).
\]
Define the unnormalised and normalised physical densities
\[
 g_i(x,s)
 :=
 h_i(x)\ind_{[0,\tau_i(x))}(s),
 \qquad
 r_i(x,s)
 :=
 \frac{g_i(x,s)}{A_i}.
\]
The probability measure with density $r_i$ is denoted
$\mu_i^{\rm ph}$.

\begin{theorem}
\label{thm:physical-exact}
For measurable roofs and densities as above,
\begin{align}
 D_{\rm ph}
 &:=
 \norm{g_1-g_0}_{L^1(\mathcal R)}
 \label{eq:Dph}
 \\
 &=
 \int_I
 \Bigl[
 \min\{\tau_0,\tau_1\}|h_1-h_0|
 +
 (\tau_1-\tau_0)_+h_1
 +
 (\tau_0-\tau_1)_+h_0
 \Bigr]
 \,\dd x.
 \nonumber
\end{align}
Consequently,
\begin{equation}\label{eq:physical-tv-normalization}
 d_{\TV}(\mu_1^{\rm ph},\mu_0^{\rm ph})
 \le
 \frac{D_{\rm ph}}{\max\{A_0,A_1\}}.
\end{equation}
Moreover,
\begin{align}
 D_{\rm ph}
 &\le
 \norm{\tau_1}_\infty
 \norm{h_1-h_0}_1
 +
 \int_Ih_0|\tau_1-\tau_0|\,\dd x,
 \label{eq:physical-bound1}
 \\
 D_{\rm ph}
 &\le
 \norm{\tau_0}_\infty
 \norm{h_1-h_0}_1
 +
 \int_Ih_1|\tau_1-\tau_0|\,\dd x,
 \label{eq:physical-bound2}
\end{align}
whenever the right-hand sides are finite.
\end{theorem}

\begin{proof}
Fix $x$ and abbreviate
\[
 a=\tau_0(x),
 \qquad
 b=\tau_1(x).
\]
On the overlap
\[
 0\le s<\min\{a,b\},
\]
the pointwise difference is
\[
 |h_1(x)-h_0(x)|.
\]
If $b>a$, the remaining strip $a\le s<b$ contributes $h_1(x)$ per unit
height; if $a>b$, the strip $b\le s<a$ contributes $h_0(x)$. Integrating
in $s$ and then in $x$ proves \cref{eq:Dph}.

Apply \cref{lem:normalization} with $f_i=g_i$ and the reference measure
$\dd x\,\dd s$ to obtain \cref{eq:physical-tv-normalization}.

For \cref{eq:physical-bound1}, write
\[
 g_1-g_0
 =
 (h_1-h_0)\ind_{[0,\tau_1)}
 +
 h_0
 \bigl(
 \ind_{[0,\tau_1)}
 -
 \ind_{[0,\tau_0)}
 \bigr).
\]
Taking the $L^1(\mathcal R)$ norm and integrating fibrewise gives
\[
 D_{\rm ph}
 \le
 \int_I\tau_1|h_1-h_0|\,\dd x
 +
 \int_Ih_0|\tau_1-\tau_0|\,\dd x,
\]
which implies \cref{eq:physical-bound1}. Interchanging the roles of
$0$ and $1$ gives \cref{eq:physical-bound2}.
\end{proof}

\begin{corollary}
\label{cor:physical-bounded}
Suppose
\[
 0<m\le\tau_i\le M<\infty
 \qquad\text{a.e.},
 \qquad i=0,1.
\]
Then
\begin{equation}\label{eq:physical-uniform}
 d_{\TV}(\mu_1^{\rm ph},\mu_0^{\rm ph})
 \le
 \frac1m
 \left(
 M\norm{h_1-h_0}_1
 +
 \min_{j\in\{0,1\}}
 \int_Ih_j|\tau_1-\tau_0|\,\dd x
 \right).
\end{equation}
In particular,
\begin{equation}\label{eq:physical-Linfty}
 d_{\TV}(\mu_1^{\rm ph},\mu_0^{\rm ph})
 \le
 \frac1m
 \left(
 M\norm{h_1-h_0}_1
 +
 \norm{\tau_1-\tau_0}_\infty
 \right).
\end{equation}
If $\norm{h_j}_\infty\le H$ for at least one $j\in\{0,1\}$, then the roof
term in \cref{eq:physical-uniform} may instead be bounded by
\[
 H\norm{\tau_1-\tau_0}_1.
\]
\end{corollary}

\begin{proof}
Since $A_i\ge m$,
\cref{eq:physical-tv-normalization,eq:physical-bound1,eq:physical-bound2}
give \cref{eq:physical-uniform}. The remaining estimates follow from
\[
 \int_Ih_j|\tau_1-\tau_0|\,\dd x
 \le
 \norm{\tau_1-\tau_0}_\infty
\]
and, when $h_j\in L^\infty$,
\[
 \int_Ih_j|\tau_1-\tau_0|\,\dd x
 \le
 H\norm{\tau_1-\tau_0}_1.
\]
\end{proof}

The dependence on roof motion is already optimal in the simplest example.

\begin{proposition}
\label{prop:constant-roof-tv}
Suppose first that the base density is
fixed and that the roof functions are constant:
\[
\tau_i\equiv c_i>0.
\]
Then
\begin{equation}\label{eq:constant-roof-tv}
 d_{\TV}(\mu_1^{\rm ph},\mu_0^{\rm ph})
 =
 \frac{|c_1-c_0|}{\max\{c_0,c_1\}}.
\end{equation}
If instead $c_0=c_1=c$ but the base densities differ, then
\begin{equation}\label{eq:common-constant}
 d_{\TV}(\mu_1^{\rm ph},\mu_0^{\rm ph})
 =
 \frac12\norm{h_1-h_0}_1.
\end{equation}
\end{proposition}

\begin{proof}
Assume $c_1\ge c_0$. The physical densities are $h/c_i$ on
$0\le s<c_i$. Thus
\begin{align*}
 2d_{\TV}
 &=
 \int_I\int_0^{c_0}
 h
 \left|
 \frac1{c_1}-\frac1{c_0}
 \right|
 \,\dd s\,\dd x
 +
 \int_I\int_{c_0}^{c_1}
 \frac{h}{c_1}
 \,\dd s\,\dd x
 \\
 &=
 2\frac{c_1-c_0}{c_1},
\end{align*}
which proves \cref{eq:constant-roof-tv}; the other order is symmetric.

If $c_0=c_1=c$, then both measures have density $h_i/c$ on the same
strip of height $c$, so
\[
 2d_{\TV}
 =
 \frac1c
 \int_I\int_0^c|h_1-h_0|\,\dd s\,\dd x
 =
 \norm{h_1-h_0}_1.
\]
\end{proof}

\subsection{Normalised-height comparison}

Let
\[
 \Mhat:=I\times[0,1),
 \qquad
 \Psi_\tau(x,u)
 =
 (x,u\tau(x)).
\]
This is a measurable coordinate identification of the half-open domains;
it is not asserted to conjugate semiflows having different roofs.

If $f$ is a density on $M_\tau$ with respect to $\dd x\,\dd s$, its
transported density on $\Mhat$ is
\begin{equation}\label{eq:Jtau}
 (J_\tau f)(x,u)
 :=
 \tau(x)f(x,u\tau(x)).
\end{equation}
Indeed,
\[
 \dd s=\tau(x)\,\dd u,
\]
so
\[
 \norm{J_\tau f}_{L^1(\Mhat)}
 =
 \norm{f}_{L^1(M_\tau)}.
\]
For the canonical suspension measure the transported density is therefore
\begin{equation}\label{eq:q}
 q(x,u)
 =
 \frac{h(x)\tau(x)}{A},
 \qquad
 A=\int_Ih\tau\,\dd x.
\end{equation}

\begin{theorem}[Stability on the fixed cylinder]
\label{thm:cylinder-stability}
Let
\[
 p_i=h_i\tau_i,
 \qquad
 A_i=\int_Ip_i\,\dd x>0.
\]
Let
\[
 \muHat_i
 =
 (\Psi_{\tau_i}^{-1})_*\mu_i
\]
be the transported suspension probability measures on $\Mhat$. Then
\begin{equation}\label{eq:cylinder-general}
 d_{\TV}(\muHat_1,\muHat_0)
 \le
 \frac{\norm{p_1-p_0}_1}{\max\{A_0,A_1\}}.
\end{equation}
Moreover,
\begin{equation}\label{eq:pbound}
 \norm{p_1-p_0}_1
 \le
 \norm{\tau_1}_\infty
 \norm{h_1-h_0}_1
 +
 \int_Ih_0|\tau_1-\tau_0|\,\dd x.
\end{equation}
If
\[
 0<m\le\tau_i\le M,
\]
then
\begin{equation}\label{eq:cylinder-uniform}
 d_{\TV}(\muHat_1,\muHat_0)
 \le
 \frac1m
 \left(
 M\norm{h_1-h_0}_1
 +
 \norm{\tau_1-\tau_0}_\infty
 \right).
\end{equation}
If $\norm{h_0}_\infty\le H$, the final roof term may be replaced by
\[
 H\norm{\tau_1-\tau_0}_1.
\]
\end{theorem}

\begin{proof}
The densities of $\muHat_i$ are $p_i/A_i$ and are independent of $u$.
Thus \cref{eq:cylinder-general} follows immediately from
\cref{lem:normalization}. The identity
\[
 p_1-p_0
 =
 \tau_1(h_1-h_0)
 +
 h_0(\tau_1-\tau_0)
\]
gives \cref{eq:pbound}. Under the common roof bounds,
\[
 A_i\ge m,
 \qquad
 \norm{\tau_1}_\infty\le M,
\]
yielding \cref{eq:cylinder-uniform}. The $L^1$ roof estimate follows by
\[
 \int_Ih_0|\tau_1-\tau_0|\,\dd x
 \le
 H\norm{\tau_1-\tau_0}_1.
\]
\end{proof}

\begin{proposition}
\label{prop:geometry-difference}
If $h$ is fixed and $\tau_i\equiv c_i>0$, then
\[
 \muHat_1
 =
 \muHat_0
 =
 h(x)\,\dd x\,\dd u,
\]
whereas \cref{eq:constant-roof-tv} gives
\[
 d_{\TV}(\mu_1^{\rm ph},\mu_0^{\rm ph})
 =
 \frac{|c_1-c_0|}{\max\{c_0,c_1\}}.
\]
Thus normalised-height total variation and physical-height total variation
are not equivalent metrics on roof perturbations.
\end{proposition}

\begin{proof}
For a constant roof $c_i$, formula \cref{eq:q} gives
\[
 q_i
 =
 \frac{hc_i}{c_i\int h}
 =
 h.
\]
The physical formula is \cref{prop:constant-roof-tv}.
\end{proof}
\begin{remark}
The two comparison procedures describe different aspects
of roof perturbations. In physical-height coordinates,
changing the roof changes the support of the suspension
measure and therefore contributes directly to its
total-variation distance. In normalised-height coordinates,
the fibre domains are identified with a fixed interval,
and a change in the roof enters only through the density
proportional to $h\tau$. Consequently, the two total-variation
distances need not be equivalent under roof perturbations.
\end{remark}
\section{Linear response under base and roof perturbations}

Let $\theta$ vary in an interval containing $0$. Suppose
\[
 h_\theta\ge0,
 \qquad
 \int_Ih_\theta\,\dd x=1,
\]
and
\begin{equation}\label{eq:hresponse}
 \frac{h_\theta-h_0}{\theta}
 \longrightarrow
 \dot h
 \qquad\text{in }L^1(I)
\end{equation}
as $\theta\to0$.

Suppose also that, for common constants $0<m<M<\infty$,
\begin{equation}\label{eq:roof-uniform-family}
 m\le\tau_\theta\le M
 \qquad\text{a.e. for }|\theta|\text{ small},
\end{equation}
and
\begin{equation}\label{eq:tauresponse}
 \frac{\tau_\theta-\tau_0}{\theta}
 \longrightarrow
 \dot\tau
 \qquad\text{in }L^\infty(I).
\end{equation}
Define
\[
 p_\theta=h_\theta\tau_\theta,
 \qquad
 A_\theta=\int_Ip_\theta\,\dd x,
 \qquad
 q_\theta=\frac{p_\theta}{A_\theta}.
\]

\begin{theorem}
\label{thm:cylinder-response}
Under
\cref{eq:hresponse,eq:roof-uniform-family,eq:tauresponse},
\begin{equation}\label{eq:pdot}
 \frac{p_\theta-p_0}{\theta}
 \longrightarrow
 \dot p
 :=
 \tau_0\dot h+h_0\dot\tau
 \qquad\text{in }L^1(I),
\end{equation}
\begin{equation}\label{eq:Adot}
 \frac{A_\theta-A_0}{\theta}
 \longrightarrow
 \dot A
 :=
 \int_I\dot p\,\dd x,
\end{equation}
and
\begin{equation}\label{eq:qdot}
 \frac{q_\theta-q_0}{\theta}
 \longrightarrow
 \dot q
 :=
 \frac{\dot p}{A_0}
 -
 \frac{p_0\dot A}{A_0^2}
 \qquad\text{in }L^1(\Mhat).
\end{equation}
In particular, $\theta\mapsto\muHat_\theta$ is differentiable in total
variation on the fixed cylinder.
\end{theorem}

\begin{proof}
Write
\begin{align*}
 \frac{p_\theta-p_0}{\theta}
 &=
 \tau_0
 \frac{h_\theta-h_0}{\theta}
 +
 h_0
 \frac{\tau_\theta-\tau_0}{\theta}
 +
 \frac{h_\theta-h_0}{\theta}
 (\tau_\theta-\tau_0).
\end{align*}
The first term converges in $L^1$ to $\tau_0\dot h$ because
$\tau_0\in L^\infty$. The second converges to $h_0\dot\tau$ in $L^1$
because the convergence in \cref{eq:tauresponse} is uniform and
$h_0\in L^1$. The $L^1$ norm of the last term is bounded by
\[
 \norm{
 \frac{h_\theta-h_0}{\theta}
 }_1
 \norm{\tau_\theta-\tau_0}_\infty,
\]
which tends to zero because the first factor remains bounded and the second
tends to zero. This proves \cref{eq:pdot}.

Integration yields \cref{eq:Adot}. Since $A_0\ge m>0$ and
$A_\theta\to A_0$, the standard quotient expansion gives
\cref{eq:qdot} in $L^1$. Independence of $u$ identifies the $L^1(I)$ and
$L^1(\Mhat)$ norms here.
\end{proof}

For physical-height response, work on the compact strip
\[
 \mathcal R_M:=I\times[0,M].
\]
For a finite signed measure $\sigma$ on $\mathcal R_M$, define
\[
 \norm{\sigma}_{\BL^*}
 :=
 \sup
 \left\{
 \left|
 \int F\,\dd\sigma
 \right|
 :
 F\in\operatorname{Lip}(\mathcal R_M),
 \ 
 \norm{F}_\infty+\operatorname{Lip}(F)\le1
 \right\}.
\]
Assume from now on that the roofs are Borel measurable, so their graphs are
Borel subsets of $\mathcal R_M$.

Define the graph-supported signed measure $\Gamma$ by
\begin{equation}\label{eq:Gamma}
 \int_{\mathcal R_M}F\,\dd\Gamma
 :=
 \int_I
 h_0(x)\dot\tau(x)
 F(x,\tau_0(x))
 \,\dd x.
\end{equation}

\begin{theorem}
\label{thm:physical-response}
Let $\mu_\theta^{\rm ph}$ be the canonical suspension measures embedded in
$\mathcal R_M$. Under
\cref{eq:hresponse,eq:roof-uniform-family,eq:tauresponse},
\begin{equation}\label{eq:physical-response-limit}
 \frac{\mu_\theta^{\rm ph}-\mu_0^{\rm ph}}{\theta}
 \longrightarrow
 \dot\mu^{\rm ph}
 \qquad\text{in }\BL^*,
\end{equation}
where
\begin{equation}\label{eq:physical-derivative}
 \dot\mu^{\rm ph}
 =
 \frac1{A_0}
 \left[
 \dot h(x)
 \ind_{[0,\tau_0(x))}(s)
 \,\dd x\,\dd s
 +
 \Gamma
 \right]
 -
 \frac{\dot A}{A_0}\mu_0^{\rm ph}.
\end{equation}

Equivalently, if
\[
 k_\theta
 :=
 \frac{h_\theta}{A_\theta},
 \qquad
 \dot k
 :=
 \frac{\dot h}{A_0}
 -
 \frac{h_0\dot A}{A_0^2},
\]
then
\begin{equation}\label{eq:physical-derivative-k}
 \dot\mu^{\rm ph}
 =
 \dot k(x)
 \ind_{[0,\tau_0(x))}(s)
 \,\dd x\,\dd s
 +
 k_0(x)\dot\tau(x)
 \,\dd x\,
 \delta_{\tau_0(x)}(\dd s).
\end{equation}
\end{theorem}

\begin{proof}
For a bounded Lipschitz test function $F$, set
\[
 B_\theta F(x)
 :=
 \int_0^{\tau_\theta(x)}
 F(x,s)\,\dd s.
\]
Then
\[
 \mu_\theta^{\rm ph}(F)
 =
 \frac1{A_\theta}
 \int_Ih_\theta B_\theta F\,\dd x.
\]

The vertical Lipschitz property gives, pointwise in $x$,
\begin{align*}
 &B_\theta F(x)-B_0F(x)
 -
 F(x,\tau_0(x))
 (\tau_\theta(x)-\tau_0(x))
 \\
 &\qquad
 =
 \int_{\tau_0(x)}^{\tau_\theta(x)}
 \bigl(
 F(x,s)-F(x,\tau_0(x))
 \bigr)\,\dd s,
\end{align*}
with the oriented-integral convention. Hence
\begin{equation}\label{eq:B-remainder}
 \norm{
 B_\theta F-B_0F-
 F(\cdot,\tau_0)(\tau_\theta-\tau_0)
 }_\infty
 \le
 \frac12
 \operatorname{Lip}(F)
 \norm{\tau_\theta-\tau_0}_\infty^2.
\end{equation}
Since \cref{eq:tauresponse} implies
\[
 \norm{\tau_\theta-\tau_0}_\infty
 =
 O(|\theta|),
\]
division by $\theta$ in \cref{eq:B-remainder} yields a remainder that is
$o(1)$ uniformly over the unit ball of $\BL$.

Now write
\begin{align*}
 &\frac1\theta
 \left(
 \int h_\theta B_\theta F
 -
 \int h_0B_0F
 \right)
 \\
 &=
 \int
 \frac{h_\theta-h_0}{\theta}
 B_0F
 +
 \int
 h_0
 \frac{B_\theta F-B_0F}{\theta}
 \\
 &\qquad
 +
 \int
 \frac{h_\theta-h_0}{\theta}
 (B_\theta F-B_0F).
\end{align*}
The first term converges uniformly over the $\BL$ unit ball to
\[
 \int\dot h\,B_0F,
\]
because
\[
 \norm{B_0F}_\infty
 \le
 M\norm{F}_\infty.
\]
By \cref{eq:B-remainder,eq:tauresponse}, the second converges uniformly to
\[
 \int_I
 h_0\dot\tau F(\cdot,\tau_0)\,\dd x,
\]
which is $\Gamma(F)$.

The last term tends uniformly to zero because
\[
 \norm{B_\theta F-B_0F}_\infty
 \le
 \norm{F}_\infty
 \norm{\tau_\theta-\tau_0}_\infty
 =
 O(|\theta|)
\]
and
\[
 \frac{h_\theta-h_0}{\theta}
\]
is bounded in $L^1$.

Finally use
\[
 A_\theta^{-1}
 =
 A_0^{-1}
 -
 \theta\frac{\dot A}{A_0^2}
 +
 o(\theta)
\]
and the uniform boundedness of the preceding functionals to obtain
\cref{eq:physical-derivative}. Rearranging the absolutely continuous terms
gives \cref{eq:physical-derivative-k}.
\end{proof}
\begin{remark}
The two terms in the first variation in \cref{eq:physical-derivative-k}
have different measure-theoretic character. The first
term is absolutely continuous with respect to planar
Lebesgue measure, whereas the second is supported on
the graph of the unperturbed roof. Since this graph
has planar Lebesgue measure zero, the second term is
singular whenever $k_0\dot{\tau}$ is nonzero on a set
of positive Lebesgue measure. This distinction is
responsible for the possible failure of differentiability
in the total-variation norm.
\end{remark}
The next result sharpens the qualitative statement that moving supports
obstruct total-variation response.

\begin{theorem}
\label{thm:metric-derivative}
Under the hypotheses of \cref{thm:physical-response}, let
\[
 k_\theta
 =
 \frac{h_\theta}{A_\theta},
 \qquad
 \dot k
 =
 \frac{\dot h}{A_0}
 -
 \frac{h_0\dot A}{A_0^2}.
\]
Then
\begin{equation}\label{eq:metric-derivative}
 \boxed{
 \lim_{\theta\to0}
 \frac{
 d_{\TV}
 (\mu_\theta^{\rm ph},\mu_0^{\rm ph})
 }{|\theta|}
 =
 \frac12
 \left(
 \int_I\tau_0(x)|\dot k(x)|\,\dd x
 +
 \int_Ik_0(x)|\dot\tau(x)|\,\dd x
 \right).
 }
\end{equation}
Thus the physical suspension-measure curve has a two-sided metric derivative
in total variation even when it has no vector-valued derivative in the
total-variation norm.
\end{theorem}

\begin{proof}
The physical densities are
\[
 r_\theta(x,s)
 =
 k_\theta(x)
 \ind_{[0,\tau_\theta(x))}(s).
\]
Repeating the fibrewise calculation used in the proof of
\cref{thm:physical-exact}, now for the nonnegative coefficients
$k_\theta$ and $k_0$, gives
\begin{align}
 \norm{r_\theta-r_0}_1
 &=
 \int_I
 \min\{\tau_\theta,\tau_0\}
 |k_\theta-k_0|
 \,\dd x
 \label{eq:metric-split}
 \\
 &\quad
 +
 \int_I
 (\tau_\theta-\tau_0)_+
 k_\theta
 \,\dd x
 +
 \int_I
 (\tau_0-\tau_\theta)_+
 k_0
 \,\dd x.
 \nonumber
\end{align}

By \cref{eq:hresponse,eq:Adot} and the quotient rule,
\[
 \frac{k_\theta-k_0}{\theta}
 \longrightarrow
 \dot k
 \qquad\text{in }L^1(I).
\]
Therefore
\[
 \frac{|k_\theta-k_0|}{|\theta|}
 \longrightarrow
 |\dot k|
 \qquad\text{in }L^1(I).
\]
Also
\[
 \min\{\tau_\theta,\tau_0\}
 \longrightarrow
 \tau_0
 \qquad\text{in }L^\infty.
\]
Hence the first term in \cref{eq:metric-split}, divided by $|\theta|$,
converges to
\[
 \int_I\tau_0|\dot k|\,\dd x.
\]

Set
\[
 a_\theta
 :=
 \frac{\tau_\theta-\tau_0}{|\theta|}.
\]
Along $\theta\downarrow0$,
\[
 a_\theta\to\dot\tau
 \qquad\text{in }L^\infty,
\]
and along $\theta\uparrow0$,
\[
 a_\theta\to-\dot\tau
 \qquad\text{in }L^\infty.
\]
In either case,
\[
 (a_\theta)_+k_\theta
 +
 (-a_\theta)_+k_0
 \longrightarrow
 k_0|\dot\tau|
 \qquad\text{in }L^1(I),
\]
because $k_\theta\to k_0$ in $L^1$ and $(a_\theta)$ is uniformly bounded
in $L^\infty$. Thus the sum of the final two terms in
\cref{eq:metric-split}, divided by $|\theta|$, converges to
\[
 \int_Ik_0|\dot\tau|\,\dd x
\]
from both sides.

Finally,
\[
 d_{\TV}(\mu_\theta^{\rm ph},\mu_0^{\rm ph})
 =
 \frac12\norm{r_\theta-r_0}_1,
\]
which proves \cref{eq:metric-derivative}.
\end{proof}
The criterion in \cref{eq:criterion} identifies the precise
role of the moving roof in total-variation differentiability.
If the roof has a nonzero first-order velocity on a set
where the unperturbed fibre density is positive, the
first variation contains a nonzero singular component.
Such a component cannot be the total-variation limit
of absolutely continuous difference quotients. If the
weighted roof velocity vanishes almost everywhere,
the moving-support contribution is negligible to first
order in $L^1$, and the derivative is given by the
absolutely continuous term in \cref{eq:tv-derivative}.

\begin{theorem}
\label{thm:tv-criterion}
Under the hypotheses of \cref{thm:physical-response}, the map
\[
 \theta\longmapsto\mu_\theta^{\rm ph}
\]
is differentiable at $0$ in the total-variation norm if and only if
\begin{equation}\label{eq:criterion}
 k_0(x)\dot\tau(x)
 =
 0
 \qquad\text{for Lebesgue-a.e. }x,
\end{equation}
or equivalently
\[
 h_0\dot\tau=0
 \qquad\text{almost everywhere}.
\]
When \cref{eq:criterion} holds, the derivative is absolutely continuous and
equals
\begin{equation}\label{eq:tv-derivative}
 \dot\mu^{\rm ph}_{\TV}
 =
 \dot k(x)
 \ind_{[0,\tau_0(x))}(s)
 \,\dd x\,\dd s.
\end{equation}
If $h_0>0$ almost everywhere, total-variation differentiability is
therefore equivalent to
\[
 \dot\tau=0
 \qquad\text{almost everywhere}.
\]
\end{theorem}

\begin{proof}
\emph{Sufficiency.}
Decompose
\begin{align*}
 r_\theta-r_0
 &=
 (k_\theta-k_0)
 \ind_{[0,\tau_0)}
 \\
 &\qquad
 +
 k_\theta
 \bigl(
 \ind_{[0,\tau_\theta)}
 -
 \ind_{[0,\tau_0)}
 \bigr).
\end{align*}
The first term divided by $\theta$ converges in
$L^1(\mathcal R_M)$ to
\[
 \dot k\ind_{[0,\tau_0)},
\]
since
\[
 \int_I
 \tau_0
 \left|
 \frac{k_\theta-k_0}{\theta}
 -
 \dot k
 \right|
 \,\dd x
 \le
 M
 \norm{
 \frac{k_\theta-k_0}{\theta}
 -
 \dot k
 }_1
 \longrightarrow0.
\]

For the moving-support term, we have
\begin{align*}
 &\frac1{|\theta|}
 \norm{
 k_\theta
 \bigl(
 \ind_{[0,\tau_\theta)}
 -
 \ind_{[0,\tau_0)}
 \bigr)
 }_1
 \\
 &=
 \int_I
 k_\theta
 \frac{
 |\tau_\theta-\tau_0|
 }{|\theta|}
 \,\dd x
 \\
 &\le
 \int_I
 k_0
 \frac{
 |\tau_\theta-\tau_0|
 }{|\theta|}
 \,\dd x
 +
 \norm{k_\theta-k_0}_1
 \norm{
 \frac{\tau_\theta-\tau_0}{\theta}
 }_\infty.
\end{align*}
By the $L^\infty$ differentiability of the roof functions,
\[
 \frac{\tau_\theta-\tau_0}{\theta}
 \longrightarrow
 \dot\tau
 \qquad\text{in }L^\infty(I).
\]
Consequently,
\[
 \frac{|\tau_\theta-\tau_0|}{|\theta|}
 \longrightarrow
 |\dot\tau|
 \qquad\text{in }L^\infty(I).
\]
It follows that
\[
 \int_I
 k_0
 \frac{|\tau_\theta-\tau_0|}{|\theta|}
 \,\dd x
 \longrightarrow
 \int_I k_0|\dot\tau|\,\dd x
 =0,
\]
where the last equality follows from \cref{eq:criterion}.
The remaining term tends to zero because
\[
 \norm{k_\theta-k_0}_1=O(|\theta|)
\]
and
\[
 \norm{
 \frac{\tau_\theta-\tau_0}{\theta}
 }_\infty
\]
is uniformly bounded for sufficiently small $|\theta|$.
Thus the moving-support contribution, divided by $|\theta|$,
converges to zero in $L^1(\mathcal R_M)$. Combining the two
terms proves total-variation differentiability, with derivative
given by \cref{eq:tv-derivative}.

\emph{Necessity.}
Every difference quotient
\[
 \frac{
 \mu_\theta^{\rm ph}
 -
 \mu_0^{\rm ph}
 }{\theta}
\]
is absolutely continuous with respect to planar Lebesgue measure
$\dd x\,\dd s$. The subspace of finite signed measures that are
absolutely continuous with respect to a fixed reference measure is
closed in total variation: if
\[
 \sigma_n\ll m
 \qquad\text{and}\qquad
 \norm{\sigma_n-\sigma}_{\mathrm{var}}\to0,
\]
then, for every measurable set $E$ with $m(E)=0$,
\[
 |\sigma(E)|
 \le
 \norm{\sigma-\sigma_n}_{\mathrm{var}}
 \to0.
\]
Therefore any total-variation derivative would have to be
absolutely continuous with respect to $\dd x\,\dd s$.

On the other hand, total-variation convergence implies
bounded-Lipschitz convergence, so the derivative must coincide
with the unique $\BL^*$ derivative in
\cref{eq:physical-derivative-k}.

The graph
\[
 G_0
 =
 \{(x,\tau_0(x)):x\in I\}
\]
is Borel and has planar Lebesgue measure zero by Fubini's theorem,
because each vertical section is a singleton. The graph term in
\cref{eq:physical-derivative-k} is therefore singular with respect
to planar Lebesgue measure. Its total variation is
\[
 \int_I
 |k_0(x)\dot\tau(x)|
 \,\dd x.
\]
It must vanish if the whole derivative is absolutely continuous.
This proves \cref{eq:criterion}.
\end{proof}

\begin{example}\label{ex:pure-roof}
Let
\[
 h_\theta=h,
 \qquad
 \tau_\theta\equiv c+\theta
\]
for $|\theta|<c/2$. Then
\[
 A_\theta=c+\theta,
 \qquad
 k_\theta=\frac{h}{c+\theta},
\]
and
\[
 d_{\TV}
 (\mu_\theta^{\rm ph},\mu_0^{\rm ph})
 =
 \frac{|\theta|}
 {\max\{c,c+\theta\}}.
\]
Hence the curve is locally Lipschitz and has metric derivative $1/c$, but
it is not differentiable in total variation. Its weak derivative is
\[
 \frac1c
 h(x)\,\dd x\,\delta_c(\dd s)
 -
 \frac1{c^2}
 h(x)
 \ind_{[0,c)}(s)
 \,\dd x\,\dd s.
\]
\end{example}

\begin{proof}
The distance formula is \cref{prop:constant-roof-tv}. The derivative
displayed above follows from \cref{eq:physical-derivative-k}, since
\[
 \dot\tau=1,
 \qquad
 \dot k=-\frac{h}{c^2}.
\]
Non-differentiability in total variation follows from
\cref{thm:tv-criterion} because
\[
 h\dot\tau=h\not\equiv0.
\]
\end{proof}

\section{Applications to piecewise-expanding maps}

The preceding theorems require no expansion hypothesis. Expansion becomes
relevant only when one asks how $h_\theta$ depends on perturbations of the
base map. 
Concrete classes of interval maps for which transfer-operator
methods establish the existence of absolutely continuous invariant
measures and quasi-compactness include piecewise-expanding
$C^{1+\varepsilon}$ maps and piecewise-convex maps with countably
many branches; see
\cite{RajputGora2025,GoraRajput2026}.
These results provide examples of base dynamics to which the
measure-theoretic suspension estimates developed here may be applied.

Consider a family
\[
 T_t:I\to I
\]
of piecewise-expanding maps with unique absolutely continuous invariant
probability measures
\[
 h_t\,\dd x.
\]
A useful abstract formulation is the following.

\begin{corollary}
\label{cor:transfer-modulus}
Assume that for $s,t$ in a parameter set,
\[
 \norm{h_t-h_s}_1
 \le
 C_h\omega_h(s,t),
\]
and that roofs satisfy common bounds
\[
 0<m\le\tau_t\le M
\]
and
\[
 \norm{\tau_t-\tau_s}_\infty
 \le
 C_\tau\omega_\tau(s,t).
\]
Then both comparison geometries satisfy
\begin{align}
 d_{\TV}(\mu_t^{\rm ph},\mu_s^{\rm ph})
 &\le
 \frac1m
 \bigl(
 MC_h\omega_h(s,t)
 +
 C_\tau\omega_\tau(s,t)
 \bigr),
 \label{eq:transfer-ph}
 \\
 d_{\TV}(\muHat_t,\muHat_s)
 &\le
 \frac1m
 \bigl(
 MC_h\omega_h(s,t)
 +
 C_\tau\omega_\tau(s,t)
 \bigr).
 \label{eq:transfer-hat}
\end{align}
If
\[
 \sup_t\norm{h_t}_\infty\le H,
\]
the roof assumption may be weakened to an $L^1$ estimate
\[
 \norm{\tau_t-\tau_s}_1
 \le
 C_\tau\omega_\tau(s,t),
\]
with $C_\tau$ replaced by $HC_\tau$ on the right-hand sides.
\end{corollary}

\begin{proof}
Apply
\cref{cor:physical-bounded,thm:cylinder-stability}.
The $L^1$ roof variant follows from the weighted estimates in those
results.
\end{proof}

Alves and Etubi \cite{AlvesEtubi2025} establish, under their
uniformity condition and additional geometric perturbation
assumptions, a modulus of continuity for the invariant
densities of families of piecewise-expanding maps. More
precisely, if the family satisfies the hypotheses of
Theorem A in \cite{AlvesEtubi2025}, then
\begin{equation}\label{eq:AE}
    \norm{h_t-h_s}_1
    \leq
    C\,\mathcal E(t-s)
    \left|\log\mathcal E(t-s)\right|.
\end{equation}
We interpret $\mathcal E(t-s)|\log\mathcal E(t-s)|$
as zero whenever $\mathcal E(t-s)=0$. Here, $\mathcal E(t-s)$ denotes the geometric perturbation
bound appearing in that theorem. Its assumptions control
the changes in branch images, the displacement of the
inverse-branch comparison maps, and the corresponding
Jacobian ratios. Thus, the estimate in \cref{eq:AE}
requires these specific geometric hypotheses and does
not follow merely from a global $C^1$ closeness of the
maps.

\begin{corollary}[A concrete piecewise-expanding suspension modulus]
\label{cor:AE}
Suppose a one-dimensional family $T_t$ satisfies the hypotheses of the
density-modulus theorem of Alves--Etubi, so that \cref{eq:AE} holds.
Suppose also
\[
 0<m\le\tau_t\le M.
\]
Then
\[
 \begin{split}
 &d_{\TV}
 (\mu_t^{\rm ph},\mu_s^{\rm ph})
 +
 d_{\TV}
 (\muHat_t,\muHat_s)
 \\
 &\qquad
 \le
 C'
 \left(
 \mathcal E(t-s)
 |\log\mathcal E(t-s)|
 +
 \norm{\tau_t-\tau_s}_\infty
 \right)
 \end{split}
\]
for a constant $C'$ independent of sufficiently close $s,t$.
If the family has a uniform $L^\infty$ bound on the invariant densities,
the final roof term may be replaced by
\[
 \norm{\tau_t-\tau_s}_1.
\]
\end{corollary}

\begin{proof}
Insert \cref{eq:AE} into
\cref{eq:transfer-ph,eq:transfer-hat}
and absorb fixed constants. The $L^1$ roof version follows from the final
statement of \cref{cor:transfer-modulus}.
\end{proof}

Linear response of the suspension measure is likewise conditional on
linear response of the base density. For piecewise-expanding interval
maps, differentiability of invariant measures is a delicate property:
in unimodal piecewise-expanding families it can depend on whether the
perturbation is tangent to the topological class; see Baladi and Smania
\cite{BaladiSmania2008,BaladiSmania2010}. The suspension formulas above do
not remove that difficulty. Rather, once an \(L^1\) derivative \(\dot h\) is available for the
base, \cref{thm:cylinder-response,thm:physical-response,thm:metric-derivative,thm:tv-criterion}
describe exactly how a simultaneously moving roof modifies the response.
\section{Invariant measures and mixing properties}

The invariant-measure calculations above are independent of whether the
suspension is mixing. A constant roof gives the simplest obstruction.

\begin{proposition}
\label{prop:constant-nonmixing}
Let
\[
 \tau\equiv c>0.
\]
Then the suspension semiflow is not mixing with respect to its canonical
suspension measure, regardless of whether the base map is mixing.
\end{proposition}

\begin{proof}
Define
\[
 F(x,s)
 =
 e^{2\pi i s/c}.
\]
This is compatible with the roof identification because
\[
 F(x,c)
 =
 1
 =
 F(Tx,0).
\]
If $n$ roof crossings occur by time $t$, then
\[
 F(\Phi^t(x,s))
 =
 e^{2\pi i(s+t-nc)/c}
 =
 e^{2\pi i t/c}F(x,s).
\]
For a constant roof, the canonical suspension measure satisfies
\[
    d\mu(x,s)=\frac{1}{c}\,d\nu(x)\,ds,
    \qquad 0\leq s<c.
\]
Consequently,
\[
    \int F\,d\mu
    =\frac{1}{c}\int_I\int_0^c
      e^{2\pi i s/c}\,ds\,d\nu(x)
    =0,
\]
and $|F|=1$ everywhere. Using the standard complex correlation
with conjugation, we obtain
\[
\begin{aligned}
C_F(t)
&=
\int (F\circ\Phi_t)\,\overline{F}\,d\mu
  -\left|\int F\,d\mu\right|^2\\
&=
e^{2\pi i t/c}\int |F|^2\,d\mu\\
&=e^{2\pi i t/c}.
\end{aligned}
\]
Thus $|C_F(t)|=1$ for every $t\geq 0$, so the suspension
semiflow is not mixing.
\end{proof}

The constant-roof example illustrates that the existence
and stability of the canonical suspension measure do not
imply mixing of the suspension semiflow. Spectral and
correlation results require additional assumptions on
the base dynamics and the roof function. For example,
Tsujii \cite{Tsujii2008} develops a spectral framework
for suspension semiflows over angle-multiplying maps,
while Daltro and Varandas \cite{DaltroVarandas2021}
establish exponential decay of correlations in specified
piecewise-expanding suspension settings. Butterley
\cite{Butterley2014} develops a functional-analytic
framework for area-expanding suspension semiflows.
These results concern dynamical and spectral properties
that are independent of the measure-comparison estimates
established here.

Similarly, Keller--Liverani perturbation theory
\cite{KellerLiverani1999} is an abstract spectral-stability mechanism; in a
concrete suspension problem one must first construct suitable strong and
weak spaces, prove uniform Lasota--Yorke-type estimates, establish weak
operator closeness under joint perturbation of the branches and roof, and
impose the dynamical assumptions needed to isolate the relevant spectral
data. None of those steps follows from total-variation stability of
invariant suspension measures.

\section{Discussion and conclusion}

The canonical suspension measure depends on the pair $(h,\tau)$
through an elementary formula, but the topology in which different
roofs are compared plays an essential role. In normalised-height
coordinates, all fibre domains are identified with a fixed cylinder,
and the corresponding density is given by
\[
\frac{h\tau}{\int_I h\tau\,dx}.
\]
This representation yields quantitative stability estimates
and an $L^1$ linear-response formula whenever the corresponding
base quantities are controlled.

In physical-height coordinates, changes in the roof move the
support of the suspension measure. This produces a graph-supported
term in the first variation, which is the obstruction to
differentiability in total variation. The exact metric-derivative
formula \cref{eq:metric-derivative} separates the contributions
of interior density variation and boundary motion, while
\cref{thm:tv-criterion} gives a necessary and sufficient condition
for total-variation differentiability in terms of the first-order
roof velocity and the unperturbed fibre density.

For piecewise-expanding bases, the stability results provide
a direct connection between perturbation estimates for invariant
densities and the corresponding suspension measures. Similarly,
available linear-response results for the base density can be
used in the suspension response formulas. The analysis of
spectral stability and decay of correlations for the suspension
semiflow remains a separate problem requiring additional
operator-theoretic and dynamical assumptions.


\begin{thebibliography}{99}

\bibitem{AlvesEtubi2025}
J.~F. Alves and O.~Etubi,
\newblock Modulus of continuity of invariant densities and entropies for
piecewise expanding maps,
\newblock \emph{Ergodic Theory Dynam. Systems}
\textbf{45} (2025), 2976--2991.

\bibitem{BaladiSmania2008}
V.~Baladi and D.~Smania,
\newblock Linear response formula for piecewise expanding unimodal maps,
\newblock \emph{Nonlinearity}
\textbf{21} (2008), 677--711.

\bibitem{BaladiSmania2010}
V.~Baladi and D.~Smania,
\newblock Alternative proofs of linear response for piecewise expanding
unimodal maps,
\newblock \emph{Ergodic Theory Dynam. Systems}
\textbf{30} (2010), 1--20.

\bibitem{Butterley2014}
O.~Butterley,
\newblock Area expanding $C^{1+\alpha}$ suspension semiflows,
\newblock \emph{Comm. Math. Phys.}
\textbf{325} (2014), 803--820.

\bibitem{DaltroVarandas2021}
D.~Daltro and P.~Varandas,
\newblock Exponential decay of correlations for Gibbs measures and
semiflows over $C^{1+\alpha}$ piecewise expanding maps,
\newblock \emph{Ann. Henri Poincar\'e}
\textbf{22} (2021), 2137--2159.

\bibitem{GalatoloLucena2020}
S.~Galatolo and R.~Lucena,
\newblock Spectral gap and quantitative statistical stability for systems
with contracting fibers and Lorenz-like maps,
\newblock \emph{Discrete Contin. Dyn. Syst.}
\textbf{40} (2020), 1309--1360.

\bibitem{Keller1985}
G.~Keller,
\newblock Generalized bounded variation and applications to piecewise
monotonic transformations,
\newblock \emph{Z. Wahrsch. Verw. Gebiete}
\textbf{69} (1985), 461--478.

\bibitem{KellerLiverani1999}
G.~Keller and C.~Liverani,
\newblock Stability of the spectrum for transfer operators,
\newblock \emph{Ann. Scuola Norm. Sup. Pisa Cl. Sci. (4)}
\textbf{28} (1999), 141--152.

\bibitem{LasotaYorke1973}
A.~Lasota and J.~A. Yorke,
\newblock On the existence of invariant measures for piecewise monotonic
transformations,
\newblock \emph{Trans. Amer. Math. Soc.}
\textbf{186} (1973), 481--488.

\bibitem{Tsujii2008}
M.~Tsujii,
\newblock Decay of correlations in suspension semi-flows of
angle-multiplying maps,
\newblock \emph{Ergodic Theory Dynam. Systems}
\textbf{28} (2008), 291--317.

\bibitem{RajputGora2025}
A.~Rajput and P.~Góra,
\newblock Quasi-compactness of Frobenius--Perron operator for
piecewise expanding $C^{1+\varepsilon}$ maps of an interval,
\newblock \emph{Stochastics and Dynamics}
\textbf{25} (2025), 2550032.
\newblock
\href{https://doi.org/10.1142/S0219493725500327}
{doi:10.1142/S0219493725500327}.

\bibitem{GoraRajput2026}
P.~Góra and A.~Rajput,
\newblock Quasi-compactness of Frobenius--Perron operator for
piecewise convex maps with countably many branches,
\newblock \emph{Dynamical Systems}
\textbf{41} (2026), no.~2, 234--251.
\newblock
\href{https://doi.org/10.1080/14689367.2026.2616234}
{doi:10.1080/14689367.2026.2616234}.

\end{thebibliography}
\end{document}